\input amstex 
\documentstyle{amsppt} 
\loadbold
\let\bk\boldkey
\let\bs\boldsymbol
\magnification=1200
\hsize=5.75truein
\vsize=8.75truein 
\hcorrection{.25truein}
\loadeusm 
\let\scr\eusm
\font\Rm=cmr17 at 17truept 

\define\Ind{\text{\rm Ind}}
\define\Hom#1#2#3{\text{\rm Hom}_{#1}({#2},{#3})} 
\define\Gal#1#2{\text{\rm Gal\hskip.5pt}(#1/#2)}
\define\sw{\text{\rm sw}} 
\define\Mid{\,\big|\,} 
\define\wt#1{\widetilde{#1}} 
\let\ge\geqslant
\let\le\leqslant

\let\vD\varDelta 
 
\let\eps\epsilon 
\let\vf\varphi
\let\vF\varPhi 
\let\vG\varGamma
 
\let\vO\varOmega

\let\vS\varSigma

\let\vT\varTheta
 
\let\vX\varXi 
\document \baselineskip=14pt \parskip=4pt plus 2pt minus 1pt 
\topmatter \nologo \nopagenumbers
\title\nofrills \Rm 
Tame multiplicity and conductor for local Galois representations
\endtitle 
\rightheadtext{Tame multiplicity and conductor} 
\author 
Colin J. Bushnell and Guy Henniart 
\endauthor 
\leftheadtext{C.J. Bushnell and G. Henniart}
\affil 
King's College London and Universit\'e de Paris-Sud 
\endaffil 
\address 
King's College London, Department of Mathematics, Strand, London WC2R 2LS, UK. 
\endaddress
\email 
colin.bushnell\@kcl.ac.uk 
\endemail
\address 
Laboratoire de Math\'ematiques d'Orsay, Univ Paris-Sud, CNRS, Universit\'e
Paris-Saclay, 91405 Orsay, France.
\endaddress 
\email 
Guy.Henniart\@math.u-psud.fr 
\endemail 
\date July 2018, revised April 2019 \enddate 
\abstract 
Let $F$ be a non-Archimedean locally compact field of residual characteristic $p$. Let $\sigma$ be an irreducible smooth representation of the absolute Weil group $\Cal W_F$ of $F$ and $\sw(\sigma)$ the Swan exponent of $\sigma$. Assume $\sw(\sigma) \ge1$. Let $\scr I_F$ be the inertia subgroup of $\Cal W_F$ and $\Cal P_F$ the wild inertia subgroup. There is an essentially unique, finite, cyclic group $\varSigma$, of order prime to $p$, so that $\sigma(\Cal I_F) = \varSigma\sigma(\Cal P_F)$. In response to a query of Mark Reeder, we show that the multiplicity in $\sigma$ of any character of $\varSigma$ is bounded by $\sw(\sigma)$. 
\endabstract 
\keywords Local field, tame multiplicity, conductor bound, primitive repre\-sentation 
\endkeywords 
\subjclassyear{2000}
\subjclass 11S37, 11S15, 22E50  \endsubjclass 
\toc 
\subhead{1.} Introduction  \endsubhead 
\subhead{2.} Group-theoretic preliminaries \endsubhead 
\subhead{3.} Conductor estimate for primitive representations \endsubhead 
\subhead{4.} Certain primitive representations \endsubhead 
\subhead{5.} An estimate of the different \endsubhead 
\subhead{6.} Proof of the main theorem \endsubhead 
\endtoc 
\endtopmatter 
\head 
1. Introduction 
\endhead 
\subhead 
1.1 
\endsubhead 
Let $F$ be a non-Archimedean, locally compact field of residual characteristic $p$. Let $\bar F$ be a separable algebraic closure of $F$ and $\scr W_F$ the Weil group of $\bar F/F$. Write $\scr I_F$ for the inertia subgroup of $\scr W_F$ and $\scr P_F$ for the wild inertia subgroup. 
\par 
Let $\sigma$ be an irreducible, smooth, complex representation of $\scr W_F$. Thus $I = \sigma(\scr I_F)$ and $P = \sigma(\scr P_F)$ are finite groups, with $P$ being the unique $p$-Sylow subgroup of $I$. The quotient $I/P$ is cyclic, of order prime to $p$. It follows readily that there is a subgroup $\vS$ of $P$ such that the quotient map $I\to I/P$ induces an isomorphism $\vS\to I/P$. Thus $\vS\cap P = 1$ and $I = \vS P$. Moreover, the subgroup $\vS$, satisfying these conditions, is uniquely determined up to conjugation by an element of $P$. (See, for instance, \cite{6} Ch\. 6 Th\. 4.1, for a full discussion.) Define the {\it tame multiplicity\/} $m(\sigma)$ of $\sigma$ by 
$$ 
m(\sigma) = \underset\chi\to{\text{\rm max}}\,\dim\Hom\vS\chi\sigma, 
$$ 
where $\chi$ ranges over the group $\widehat\vS$ of linear characters of $\vS$. The integer $m(\sigma)$ does not depend on the choice of $\vS$ and, in all cases, $m(\sigma) \ge1$. 
\par 
Let $\sw(\sigma)$ be the Swan exponent of $\sigma$. We prove:
\proclaim{Tame Multiplicity Theorem} 
Let $\sigma$ be an irreducible smooth representation of $\scr W_F$. If $\sw(\sigma) > 0$, then 
$$ 
m(\sigma)\le \sw(\sigma). 
\tag 1.1.1 
$$ 
In particular, the space of $\vS$-fixed points in $\sigma$ has dimension at most $\sw(\sigma)$. 
\endproclaim 
In his paper \cite{10}, Mark Reeder gives compelling reasons for being interested in the invariant $m(\sigma)$ and the inequality (1.1.1). He proves the theorem when $\sigma$ is either essentially tame or of epipelagic type, in the sense that $\sw(\sigma) = 1$.  This paper is written in response to his query as to whether it might hold in general. 
\remark{Remarks} 
A couple of cases can be dispatched straightaway. 
\roster 
\item 
If $\sw(\sigma) = 0$, then $\vS = I$ and $\sigma$ is induced from a tamely ramified character of $\scr W_K$, where $K/F$ is unramified. It follows that $m(\sigma) = 1$. 
\item 
If $\dim\sigma = 1$ and $\sw(\sigma) \ge 1$, then $m(\sigma) = 1$ for trivial reasons. 
\endroster 
\endremark 
\subhead 
1.2 
\endsubhead 
This is, obviously, a ``small conductor'' problem: certainly $m(\sigma) \le \dim\sigma$ while, for the vast majority of representations $\sigma$, one has $\dim\sigma < \sw(\sigma)$. On the other hand, if $\sw(\sigma) = 1$ then $m(\sigma) = 1$ \cite{2}, \cite{10}. It is the contrast between these two extremes that dictates the flavour of the paper. In many cases, rather coarse estimates should suffice to give the result but, in others, delicacy is likely to be required. 
\par 
The small-conductor aspect suggests that primitive representations $\sigma$ must play a central r\^ole. At first glance, one might hope to prove the theorem for primitive representations and then proceed by induction. That light-hearted approach falls at the first hurdle. If one tries to calculate $m(\sigma)$ directly from the description of $\sigma$ in Koch's structure theory \cite{9}, the combinatorics rapidly get out of hand. Further, we have an uncertain grasp of the relation between Koch's description of $\sigma$ and the value of $\sw(\sigma)$. Examples show that there is sometimes no room for any sloppiness in the estimates. 
\par 
More positively, there is a strong lower bound for $\sw(\sigma)$ in \cite{7}. On the other side, help comes from a rather different source. Glauberman's general theory of character correspondences for finite groups \cite{5}, as developed in Isaacs' book \cite{8}, leads to an exact and manageable formula for $m(\sigma)$, but only for a restricted class of primitive representations $\sigma$. To outline this, we need some terminology. 
\par 
Let $\tau$ be an irreducible representation of $\scr W_F$. Say that $\tau$ is {\it absolutely ramified\/} if it factors through $\Gal EF$, where $E/F$ is a finite, totally ramified field extension. Let $\sigma$ be primitive and absolutely ramified, viewed as a faithful representation of $G = \Gal EF$. Let $\Gal EK$ be the centre of $G$ and let $T/F$ be the maximal tame sub-extension of $E/F$. We may reduce to the case where $\Gal ET$ is a $p$-group, and therefore the wild inertia subgroup of $G$. Let $\vS$ be a complement of $\Gal ET$ in $G$. For the purposes of this introduction, say that $\sigma$ is {\it $\vS$-homogeneous\/} if the $G$-centralizer of any non-trivial element of $\vS$ is $\vS.\Gal EK$. If $\sigma$ is $\vS$-homogeneous, then \cite{8} gives an exact formula for the character $\roman{tr}\,\sigma\Mid\vS$. 
\par 
If $\sigma$ is absolutely ramified and $\vS$-homogeneous, comparison of the character formula with the conductor bound of \cite{7} yields the theorem. This is a case in which $m(\sigma)$ can be close to $\sw(\sigma)$ (see 4.5). More generally, an absolutely ramified primitive representation is essentially a tensor product of homogeneous ones. A relatively relaxed estimate then gives the theorem in this case. 
\par 
For the third step, we prove the theorem when $\sigma$ is absolutely ramified (but not necessarily primitive). We can assume that $\sigma$ is an induced representation $\Ind_{K/F}\,\tau$, where $\tau$ is an absolutely ramified representation of $\scr W_K$ with $K \neq F$ and $m(\tau) \le \sw(\tau)$. A standard property asserts that $\sw(\sigma) = \sw(\tau)+ w_{K/F}\dim\tau$, where $w_{K/F}$ is the {\it wild exponent\/} of the extension $K/F$. The relation between $m(\sigma)$ and $m(\tau)$ is group-theoretic in nature, so we have to estimate the arithmetic quantity $w_{K/F}$ in group-theoretic terms. A rather coarse argument suffices. It shows that, relative to induction of representations, $\sw(\sigma)$ grows much more quickly that $m(\sigma)$ and so justifies the initial emphasis on primitive representations. From there on, the general case of the theorem follows easily. 
\subhead 
1.3 
\endsubhead 
The paper is arranged as follows. The necessary material from finite group theory is assembled in section 2. In section 3, we review some properties of primitive representations leading to the conductor estimate of \cite{7} Th\'eor\`eme 1.8. We give a complete proof of that result. It uses the same ideas as \cite{7} but, in the present limited context, they can be expressed more succinctly and transparently. Section 4 is the heart of the argument, proving the theorem for ``$\vS$-homogeneous'', absolutely ramified, primitive representations, as sketched above. Section 5 is the group-theoretic estimate of the wild exponent, and section 6 completes the proof. 
\remark{Acknowledgement} 
We thank the referee for his detailed comments on an earlier version. These led us to produce a much improved version. In particular, he noticed the elementary error corrected in 2.1 below. 
\endremark 
\head 
2. Group-theoretic preliminaries 
\endhead 
We gather some techniques from the representation theory of finite groups. This section has its own scheme of notation. 
\subhead 
2.1 
\endsubhead 
We consider a special class of finite $p$-groups, using the terminology of \cite{4}.  
\definition{Definition} 
Let $P$ be a finite $p$-group with centre $Z \neq P$. It is called {\it H-cyclic\/} if it satisfies the following conditions. 
\roster 
\item The centre $Z$ is cyclic, and 
\item the quotient $V = P/Z$ is elementary abelian. 
\endroster 
\enddefinition 
For convenience, we summarize the main properties of these groups, following the account in \cite{4}. For $x,y\in P$, the commutator $[x,y]$ lies in the centre $Z$ and satisfies $[x,y]^p = [x^p,y] = 1$. 
\par 
We think of the elementary abelian $p$-group $V$ as a vector space over the field $\Bbb F_p$ of $p$ elements. Let $\zeta$ be a faithful character of $Z$. The commutator pairing 
$$ 
(x,y) \longmapsto \zeta(xyx^{-1}y^{-1}), \quad x,y \in P, 
\tag 2.1.1 
$$ 
takes its values in the group $\bs\mu_p(\Bbb C)$ of complex $p$-th roots of unity. Composing with a fixed isomorphism $\bs\mu_p(\Bbb C) \to \Bbb F_p$, the pairing (2.1.1) induces an alternating bilinear form 
$$ 
h_\zeta: V\times V \longrightarrow \Bbb F_p. 
\tag 2.1.2 
$$ 
Because $Z$ is the centre of $P$, the form $h_\zeta$ is nondegenerate. Consequently, $V$ has $p^{2r}$ elements, for an integer $r\ge1$. 
\par 
A subspace $W$ of $V$ is a {\it Lagrangian subspace\/} of the alternating space $(V,h_\zeta)$ if it has exactly $p^r$ elements and $h_\zeta(w_1,w_2) = 0$ for all $w_1,w_2\in W$. 
\proclaim{Lemma} 
There is a unique irreducible representation $\tau$ of $P$ such that $\tau\Mid Z$ contains $\zeta$. It has the following properties. 
\roster 
\item The representation $\tau$ is faithful, it satisfies $\dim\tau = p^r$, and $\tau\Mid Z$ is a multiple of $\zeta$. 
\item 
Let $W$ be a Lagrangian subspace of $(V,h_\zeta)$ with inverse image $\wt W$ in $P$. The group $\wt W$ is abelian and the character $\zeta$ of $Z$ admits extension to a character $\zeta^W$ of\/ $\wt W$. For any such $\zeta^W$, one has 
$$ 
\tau \cong \Ind_{\wt W}^P\,\zeta^W. 
$$ 
\endroster 
\endproclaim 
\demo{Proof} 
See \cite{4} 8.1 Proposition. \qed 
\enddemo 
\remark{Remark} 
A finite $p$-group is called {\it extra special of class\/} $2$ if it is H-cyclic and its centre has order $p$ ({\it cf\.} \cite{6} p\. 183). More generally, let $P$ be H-cyclic with centre $Z$. Since the representation $\tau$ of the lemma is faithful, one may identify $P$ with $\tau(P)$. One can then follow Rigby's argument in \cite{11} Theorem 2 to show that $P$ is the central product of the finite cyclic $p$-group $Z$ and an extra special $p$-group of class $2$. 
\endremark 
\remark{Corrigendum} 
In the preamble to \cite{4} 8.1, we assert that an H-cyclic group is extra special of class $2$. The arguments of \cite{4} section 8 are conducted axiomatically, so this error has no effect on the results or their proofs. In particular, the lemma above remains valid. 
\endremark 
\subhead 
2.2 
\endsubhead 
Let $P$ be a finite, H-cyclic $p$-group with centre $Z$, write $V = P/Z$ and let $|V| = p^{2r}$. Let $\zeta$ be a faithful character of $Z$. We introduce another element of structure. 
\definition{Definition} 
Let $S$ be a cyclic group of automorphisms of $P$, such that 
\roster 
\item the order $|S|$ of $S$ is not divisible by $p$ and 
\item $S$ acts trivially on $Z$. 
\endroster 
\enddefinition 
Because of condition (2), the action of $S$ on $P$ fixes the commutator pairing (2.1.2), so the $\Bbb F_p$-repre\-sentation of $S$ provided by $V$ is {\it symplectic.} We consider a specific representation of the semi-direct product $G = S\ltimes P$. 
\proclaim{Lemma} 
Let $G = S \ltimes P$ and let $\tau$ be the unique irreducible representation of $G$ such that $\tau\Mid Z$ contains $\zeta$. 
\roster 
\item 
There exists a unique representation $\tilde\tau$ of $G$ such that $\tilde\tau\Mid P \cong \tau$ and $\det \tilde\tau(s) = 1$, for all $s\in S$. 
\item 
An irreducible representation $\rho$ of $G$ satisfies $\rho\Mid P \cong \tau$ if and only if there is a character $\chi$ of $S = G/P$ such that $\rho \cong \chi\otimes\tilde\tau$. 
\endroster 
\endproclaim 
\demo{Proof} 
This is a pleasant exercise, written out in \cite{1} (8.4.1) Proposition. For a very general result of this kind, see \cite{8} 13.3 Lemma. \qed 
\enddemo 
Under certain circumstances, one can write down the character $\roman{tr}\,\tilde\tau$ of $\tilde\tau$ on elements of $S$. 
\proclaim{Proposition} 
Suppose that, for every $s\in S$, $s\neq1$, the $G$-centralizer of $s$ is $SZ$. There is then a constant $\eps = \pm1$ and a character $\mu$ of $S$, such that $\mu^2 = 1$, with the following properties. 
\roster 
\item 
$\roman{tr}\,\tilde\tau(sz) = \eps\mu(s)\zeta(z)$, for $s\in S$, $s\neq 1$, and $z\in Z$. 
\item 
$p^r{-}\eps = k|S|$, for an integer $k$. 
\item 
$\roman{tr}\,\tilde\tau\Mid S = k\,\roman R_S +\eps\mu$, where $\roman R_S$ is the character of the regular representation of $S$. 
\item 
The character $\mu$ is non-trivial if and only if $|S|$ is even and $k$ is odd. 
\endroster 
\endproclaim 
\demo{Proof} 
This is a special case of \cite{8} Theorem 13.32: to translate the notation, our $\tau$ is $\chi$ in \cite{8}, while $\tilde\tau$ is $\hat\chi$ and $\zeta$ is $\beta$. Otherwise, conventions are the same. \qed 
\enddemo  
\remark{Remarks} 
\roster 
\item 
The formulas in the proposition show that the character $\roman{tr}\,\tilde\tau\Mid S$ of $S$ is determined by the group orders $|S|$ and $|V|$. Indeed, if $|S| \ge3$, the invariants $k$, $\eps$, $\mu$ are individually determined by the group orders. When $|S| = 2$, the character is determined but the invariants are not. In all cases, the character $\roman{tr}\, \tilde\tau \Mid S$ depends only on the {\it linear\/} $\Bbb F_p$-representation of $S$ afforded by $V$. 
\item 
Let $S'$ be a cyclic group, of order prime to $p$, equipped with a surjective homomorphism $S'\to S$. One may inflate $\tilde\tau$ to a representation $\tilde\tau'$ of $S'\ltimes P$ and then use the proposition to write down $\roman{tr}\,\tilde\tau'\Mid S'$. 
\endroster 
\endremark 
Character computations of this sort feature in \cite{1}, especially (8.6.1) Theorem, and have been widely used. However, the account in \cite{1} deals only with the case where the symplectic $\Bbb F_pS$-representation $P/Z$ is indecomposable. The proposition gives the exact formula for a wider class of cases. It neatly avoids an estimation process at a point where absolute precision is essential ({\it cf\.} 4.4, 4.5 below). 
\subhead 
2.3 
\endsubhead 
Suppose that the space $V$ of 2.1 decomposes as a direct sum of nonzero subspaces $V_1$, $V_2$, orthogonal with respect to the alternating form (2.1.2), say 
$$ 
V = V_1\perp V_2. 
\tag 2.3.1 
$$ 
Let $P_i$ be the inverse image of $V_i$ in $P$. The commutator group $[P_1,P_2]$ is trivial, that is, $P_1$ commutes with $P_2$. Moreover, each $P_i$ is H-cyclic with centre $Z$. 
\par 
The obvious map $P_1\times P_2 \to P$ is a surjective homomorphism with kernel $\{(z,z^{-1}): z\in Z\}$. That is, $P$ is the central product of its subgroups $P_1$, $P_2$. As in 2.1 Lemma, the group $P_i$ admits a unique irreducible representation $\tau_i$ containing the character $\zeta$ of $Z$. The representation $\tau_1\otimes \tau_2$ factors through the quotient map $\pi:P_1\times P_2 \to P$ and so $\tau\circ\pi \cong  \tau_1\otimes \tau_2$: one may reasonably write 
$$ 
\tau = \tau_1\otimes \tau_2. 
\tag 2.3.2 
$$ 
\subhead 
2.4 
\endsubhead 
Let $S$ be a cyclic group of automorphisms of $P$, as in 2.2 Definition, and suppose that the factors $V_i$ in (2.3.1) are $S$-invariant. It follows that the subgroups $P_i$ of $P$ are normalized by $S$. Let $S_i$ be the image of $S$ in $\roman{Aut}\,P_i$. Following the procedure of 2.2 Lemma, we form the representation $\tilde\tau_i$ of $S_i\ltimes P_i$. We inflate $\tilde\tau_i$ to a representation $\tilde\tau^S_i$ of $S\ltimes P_i$. We can equally set $\tau = \tau_1\otimes\tau_2$ as in (2.3.2) and extend it to a representation $\tilde\tau$ of $S\ltimes P$ as before. We then have 
$$ 
\roman{tr}\,\tilde\tau(s) = \roman{tr}\,\tilde\tau_1^S(s)\cdot \roman{tr}\,\tilde\tau_2^S(s), \quad s\in S. 
\tag 2.4.1 
$$ 
\head 
3. Conductor estimate for primitive representations
\endhead 
We give a lower bound, in terms of ramification structure, for the Swan exponent of a certain class of representations of the Weil group. Before we start, we lay down some notation and conventions to remain in force for the rest of the paper. 
\example{Notation and conventions} 
\roster 
\item 
Let $\scr W_F$ be the Weil group of a chosen separable closure $\bar F/F$. When speaking of a ``representation of $\scr W_F$'' we mean a ``smooth complex representation of $\scr W_F$''. Let $\scr I_F$ be the inertia subgroup of $\scr W_F$ and $\scr P_F$ the wild inertia subgroup. 
\item 
Let $\frak p_F$ be the maximal ideal of the discrete valuation ring in $F$. If $k\ge1$ is an integer, then $U^k_F$ is the unit group $1{+}\frak p_F^k$. The residue field of $F$ is $\Bbbk_F$. 
\item 
We use the conventions of \cite{12} when dealing with ramification groups, their numberings and the Herbrand functions $\vf$, $\psi$. 
\endroster 
\endexample
\subhead 
3.1 
\endsubhead 
An irreducible representation $\sigma$ of $\scr W_F$ is called {\it primitive\/} if $\dim\sigma > 1$ and if $\sigma$ is not induced from a representation of $\scr W_K$, where $K/F$ is a finite field extension with $K \neq F$. 
\proclaim{Hypothesis}  
For the rest of this section, we suppose that the representation $\sigma$ is primitive. 
\endproclaim 
The restriction $\sigma\Mid \scr P_F$ is then irreducible and the finite $p$-group $\sigma(\scr P_F)$ is H-cyclic \cite{11}, Theorem 1. Consequently, $\dim\sigma = p^r$,  for some $r\ge1$. Let $\bar\sigma$ be the projective representation defined by $\sigma$ and set $\scr W_K = \roman{Ker}\,\bar\sigma$. In particular, $\sigma(\scr W_K)$ is the centre of $\sigma(\scr W_F)$, so $\sigma\Mid \scr W_K$ is a multiple of a character $\zeta_\sigma$ of $\scr W_K$. 
\par 
Let $T/F$ be the maximal tamely ramified sub-extension of $K/F$. The group $\vD = \Gal KT$ is elementary abelian of order $p^{2r}$. Since $\sigma(\scr P_F)$ is H-cyclic, the pairing 
$$ 
(x,y) \longmapsto \zeta_\sigma(xyx^{-1}y^{-1}), \quad x,y\in \scr W_T, 
$$ 
induces a bilinear form 
$$ 
h_\sigma: \vD\times \vD \longrightarrow \Bbb F_p 
\tag 3.1.1 
$$ 
that is alternating and nondegenerate. The natural action of $\vT = \Gal TF$ on $\vD$ fixes $h_\sigma$, so $(\vD,h_\sigma)$ provides a {\it symplectic\/} representation of $\vT$ over the field $\Bbb F_p$. A crucial point is the following \cite{9} Theorem 4.1. 
\proclaim{Proposition} 
The symplectic $\Bbb F_p$-representation of $\vT$ on $\vD$ is \rom{$\vT$-anisotropic,} in that $\vD$ has no non-zero $\vT$-subspace on which $h_\sigma$ is identically zero. 
\endproclaim  
It is usually convenient to impose a further normalization. 
\proclaim{Lemma} 
There is a tamely ramified character $\chi$ of $\scr W_F$, such that the representation $\sigma' = \chi\otimes \sigma$ has the following properties. 
\roster 
\item The kernel of $\sigma'$ is of the form $\scr W_E$, where $E/K$ is finite, cyclic and totally wildly ramified. 
\item The order of the character $\zeta_{\sigma'}$ is finite and a power of $p$, with $\scr W_E = \roman{Ker}\,\zeta_{\sigma'}$. 
\endroster 
\endproclaim 
\demo{Proof} 
We construct the character $\chi$ in stages. First, there is an unramified character $\chi_1$ of $\scr W_F$ such that the representation $\sigma_1 = \chi_1\otimes \sigma$ has finite image. The character $\det\sigma_1$ therefore has finite order. There exists a character $\chi_2$ of $\scr W_F$, of finite order relatively prime to $p$, such that $\chi_2^{p^r}\det\sigma_1$ has finite $p$-power order. In particular, $\chi_2$ is tamely ramified. Set $\sigma_2 = \chi_2\otimes \sigma_1$, so that $\det\sigma_2$ has finite $p$-power order. The restriction of $\sigma_2$ to $\scr W_K$ is a multiple of the character $\zeta_2 = \zeta_\sigma \cdot \chi_2\chi_1\Mid \scr W_K$. By construction, $\zeta_2$ has finite $p$-power order. 
\par 
Let $\roman{Ker}\,\zeta_2 = \scr W_{E_2}$. Thus $E_2/K$ is a finite, cyclic $p$-extension. Viewing $\zeta_2$ as a character of $K^\times$ via class field theory, the extension $E_2/K$ is totally ramified if and only if $\zeta_2(K^\times) = \zeta_2(U_K)$ or, equivalently, there is a Frobenius element $\phi$ of $\scr W_K$ such that $\zeta_2(\phi) = 1$. So, suppose we have a Frobenius $\phi$ for which $\zeta_2(\phi) \neq1$. There is an unramified character $\psi$ of $\scr W_K$, of finite, $p$-power order, such that $\psi\zeta_2(\phi) =1$. This character $\psi$ is the restriction of an unramified character $\chi_3$ of $\scr W_F$ of finite, $p$-power order. Write $\zeta_3 = \chi_3\zeta_2$ and $\scr W_E = \roman{Ker}\,\zeta_3$. The extension $E/K$ is cyclic and totally ramified of $p$-power degree. Moreover, $\scr W_E = \roman{Ker}\,\sigma_3$, where $\sigma_3 = \chi_3\otimes\sigma_2$, and all assertions have been proved for $\sigma' = \sigma_3$. \qed 
\enddemo 
\remark{Remark} 
Replacing $\sigma$ by $\sigma'$ has no effect on the pairing $h_\sigma$ or the fields $K$, $T$. The Tame Multiplicity Theorem holds for $\sigma$ if and only if it holds for $\sigma'$. 
\endremark 
\subhead 
3.2 
\endsubhead 
In the notation of 3.1, we analyze the symplectic $\Bbb F_p$-representation of $\vT$ provided by $\vD$. Let $J_{K/T}$ be the set of ramification jumps of $K/T$, in the upper numbering. Since $K/T$ is abelian, these jumps are positive integers, by the Hasse-Arf Theorem \cite{12}  V Th\'eor\`eme 1. Observe that, for a real number $x\ge0$, the ramification group $\vD^x$ is an $\Bbb F_p\vT$-subspace of $\vD$. 
\proclaim{Proposition} 
Let $j\in J_{K/T}$. 
\roster 
\item 
The restriction of $h_\sigma$ to $\vD^j$ is nondegenerate. 
\item 
If\/ $W^j$ denotes the $h_\sigma$-orthogonal complement of $\vD^{1+j}$ in $\vD^j$, then $\vD$ is the orthogonal sum of the spaces $W^j$, $j\in J_{K/T}$. 
\endroster 
\endproclaim 
\demo{Proof} 
If $X$ is a subspace of $\vD$,  let $X^\perp$ be its $h_\sigma$-orthogonal complement in $\vD$. For an integer $j\ge1$, the radical of the alternating form $h_\sigma\Mid \vD^j\times\vD^j$ is $\vD^j\cap (\vD^j)^\perp$. This is an $\Bbb F_p\vT$-subspace of $\vD$ on which $h_\sigma$ is null. Since $h_\sigma$ is $\vT$-anisotropic, $\vD^j\cap (\vD^j)^\perp = 0$ whence (1) follows. 
\par 
If $j\in J_{K/T}$, then $\vD^{1+j}$ is trivial or equal to $\vD^{j'}$, where $j'$ is the least element of $J_{K/T}$ strictly greater than $j$. Assertion (2) now follows from (1). \qed 
\enddemo 
\subhead 
3.3 
\endsubhead 
We continue with the notation of 3.1, 3.2 to establish a lower bound on the Swan exponent $\sw(\sigma)$. 
\par 
First, we specify a family of Lagrangian subspaces of the alternating space $(\vD,h_\sigma)$. For each $j\in J_{K/T}$, let $\vX(j)$ be a Lagrangian subspace of the nondegenerate space $W^j$. The various $\vX(j)$ are mutually orthogonal, and so $\vX = \sum_j\vX(j)$ is Lagrangian. A Lagrangian subspace of this form will be called {\it $J$-split.} 
\proclaim{Theorem} 
Let $\vX$ be a $J$-split Lagrangian subspace of $\vD$. If $K^\vX = L$, then $J_{L/T} = J_{K/T}$. If $j_\infty$ is the largest element of $J_{L/T}$ and $e(T|F) = e$ then 
$$ 
e\,\sw(\sigma) \ge \psi_{L/T}(j_\infty) + p^rj_\infty \ge (1{+}p^r)j_\infty. 
\tag 3.3.1 
$$ 
\endproclaim 
\demo{Proof} 
We may assume, without loss, that the representation $\sigma$ has been normalized as in 3.1 Lemma. In particular, $\roman{Ker}\,\sigma = \scr W_E$, where $E/K$ is cyclic and totally wildly ramified. The extension $E/F$ is Galois. 
\par 
By construction, the extensions $K/T$ and $L/T$ have the same jumps, $J_{L/T} = J_{K/T}$. Let $\wt\vD = \Gal ET$, $\wt\vX = \Gal EL$. Since $\vX$ is a Lagrangian subspace of $\vD$, the extension $E/L$ is abelian and totally wildly ramified. The Artin Reciprocity isomorphism therefore induces a surjective homomorphism 
$$ 
a_L: U^1_L \longrightarrow \wt\vX = \Gal EL. 
$$ 
\indent 
Let $x \in \vD^{j_\infty}\cap \vX$ be non-trivial, and choose $y\in \vD^{j_\infty}$ such that $h_\sigma(x,y) \neq 0$. We have $\vD^{j_\infty} = \vD_{k_\infty}$, where $k_\infty = \psi_{K/T}(j_\infty)$, and so $x$ is an element of $\vD_{k_\infty}\cap \vX = \vX_{k_\infty}$. However, 
$$ 
\vX_{k_\infty} = \vX^{\vf_{K/L}(k_\infty)} = \vX^{\psi_{L/T}(j_\infty)}, 
$$ 
as follows from the transitivity relation $\psi_{K/T} = \psi_{K/L}\circ \psi_{L/T}$. Choose an inverse image $\tilde x$ of $x$ in $\wt \vX^{\psi_{L/T}(j_\infty)}$. As Galois operator on $E$ therefore, we have $\tilde x = a_L(v)$, for some $\psi_{L/T}(j_\infty)$-unit $v$ of $L$ (by the higher ramification theorem of local class field theory \cite{12} XV Th\'eor\`eme 1 Corollaire 3). 
\par 
On the other hand, $y$ acts on $L$ as an element of 
$$ 
(\wt\vD/\wt\vX)^{j_\infty} = (\vD/\vX)^{j_\infty} = (\vD/\vX)_{\psi_{L/T}(j_\infty)}. 
$$ 
The definition of the lower ramification sequence implies that, if $z\in U^k_L$, for some $k\ge1$, then $z^y/z$ is a $(k{+}\psi_{L/T}(j_\infty))$-unit of $L$. 
\par 
Choose an inverse image $\tilde y$ of $y$ in $\wt\vD^{j_\infty}$. Therefore 
$$ 
\tilde y^{-1}\tilde x \tilde y\tilde x^{-1} = a_L(v^{\tilde y}v^{-1}) = a_L(u), 
$$ 
where $u = v^{\tilde y}v^{-1} = v^yv^{-1}$ is a $2\psi_{L/T}(j_\infty)$-unit of $L$. 
\par 
Set $\sigma\Mid \scr W_T = \tau$. The representation $\tau$ is irreducible. Since $\vX$ is Lagrangian, $\tau$ is induced from a character $\phi$ of $\scr W_L$ extending the character $\zeta_\sigma$ of $\scr W_K$ (2.1 Lemma).  By construction, $\zeta_\sigma[y^{-1},x] = \zeta_\sigma[\tilde y^{-1},\tilde x] \neq 1$. So, if we view $\phi$ as a character of $L^\times$ via class field theory, it is non-trivial on $2\psi_{L/T}(j_\infty)$-units of $L$. That is, $\sw(\phi) \ge 2\psi_{L/T}(j_\infty)$. Let $w_{L/T}$ be the wild exponent of the extension $L/T$ ({\it cf\.} (5.1.1) below). The standard induction formula reads 
$$ 
\align 
\sw(\tau) &= \sw(\phi) + w_{L/T} \\ 
&\ge 2\psi_{L/T}(j_\infty) + w_{L/T}. 
\endalign 
$$ 
Since $[L{:}T] = p^r$ and $j_\infty$ is the largest jump of $L/T$, we have 
$$ 
\psi_{L/T}(j_\infty) = p^rj_\infty - w_{L/T} 
$$ 
by \cite{4} 1.6 Proposition. It follows that $\sw(\tau) \ge \psi_{L/T}(j_\infty) + p^rj_\infty$. The Herbrand function satisfies $\psi_{L/T}(x) \ge x$, for all $x\ge0$, so we further have 
$$ 
\sw(\tau) \ge \psi_{L/T}(j_\infty) + p^rj_\infty \ge (1{+}p^r)j_\infty. 
$$ 
Since $\sw(\tau) = e\,\sw(\sigma)$, we are done. \qed 
\enddemo 
\subhead 
3.4 
\endsubhead 
The theorem of 3.3, and its proof, apply unchanged in greater generality. We shall not use the fact here, but this is a convenient place to record it. Suppose only that the irreducible representation $\sigma$ is H-cyclic, in the sense of \cite{4}: this means that $\sigma\Mid \scr P_F$ is irreducible and that the finite $p$-group $\sigma(\scr P_F)$ is H-cyclic in the sense of 2.1. We can use all the same notation relative to $\sigma$. The inequalities (3.3.1) then hold, {\it provided the alternating form $h_\sigma$ is nondegenerate on $\vD^{j_\infty}$.} 
\head 
4. Certain primitive representations 
\endhead 
In this section, we prove the Tame Multiplicity Theorem for a certain class of primitive representations of $\scr W_F$. 
\subhead 
4.1 
\endsubhead 
Let $\sigma$ be a primitive irreducible representation of $\scr W_F$. Say that $\sigma$ is called {\it absolutely ramified\/} if the associated projective representation $\bar\sigma$ factors through a finite Galois group $\Gal LF$ for which $L/F$ is totally ramified. 
\proclaim{Theorem} 
If $\sigma$ is an irreducible, primitive, absolutely ramified representation of $\scr W_F$, then $m(\sigma) \le \sw(\sigma)$. 
\endproclaim 
The proof will occupy the rest of the section. 
\subhead 
4.2  
\endsubhead 
We normalize $\sigma$ as permitted by (3.1) Lemma and use the notation devel\-oped in 3.1. Thus $\roman{Ker}\,\bar \sigma = \scr W_K$, where $K/F$ is totally ramified. Let $T/F$ be the maximal tame sub-extension of $K/F$. In addition, $\roman{Ker}\,\sigma = \scr W_E$ where $E/K$ is cyclic and totally wildly ramified. 
\par 
Set $\vG = \Gal KF$, $\vD = \Gal KT$ and $\vT = \Gal TF$. Therefore $\vD$ is elementary abelian of order $p^{2r} = (\dim\sigma)^2$ and $\vT$ is cyclic of order prime to $p$. The restriction of $\sigma$ to $\scr W_K$ is a multiple of a character $\zeta_\sigma$ and the group $\wt\vD = \Gal ET$ is an H-cyclic $p$-group with centre $\Gal EK$. The subgroup $\vD$ admits a complement $\vS$ in $\vG$. Thus $\vS\cap \vD = \{1\}$ and $\vG = \vS\vD$. Restriction of operators induces an isomorphism $\vS \cong \vT$. In particular, $\vS$ is cyclic of order $e = e(T|F)$. 
\par 
Let $h_\sigma$ be the commutator pairing as in (3.1.1). The pair $(\vD, h_\sigma)$ affords an anisotropic, symplectic $\Bbb F_p$-representation of $\vS$, of dimension $2r$. We review the classification of such representations, following \cite{1}. 
\par 
Choose an algebraic closure $\bar{\Bbb F}_p/\Bbb F_p$, and write $\Gal{\bar{\Bbb F}_p}{\Bbb F_p} = \vO$. Let $\chi:\vS \to \bar{\Bbb F}_p^\times$ be a homomorphism and let $\Bbb F_p(\chi)$ be the field generated by the values $\chi(s)$, $s\in \vS$. The group $\vS$ acts on $\Bbb F_p(\chi)$ via the character $\chi$, that is, 
$$ 
s:x\longmapsto \chi(s)x,\quad s\in \vS,\ x\in \Bbb F_p(\chi). 
$$ 
The group $\vO$ acts on $\Hom{}\vS{\bar{\Bbb F}_p}$ in a natural way. The map $\chi\mapsto \Bbb F_p(\chi)$ then induces a bijection between $\vO\backslash \Hom{}\vS{\bar{\Bbb F}_p}$ and the set of isomorphism classes of irreducible $\Bbb F_p$-representations of $\vS$. 
\proclaim{Proposition} 
\roster 
\item 
For $\chi\in \Hom{}\vS{\bar{\Bbb F}_p^\times}$, the following conditions are equivalent. 
\itemitem{\rm (a)} The representation $\Bbb F_p(\chi)$ is symplectic, that is, it admits a nondegenerate, $\vS$-invariant, alternating form. 
\itemitem{\rm (b)} The character $\chi^{-1}$ is $\vO$-conjugate, but not equal, to $\chi$. 
\itemitem{\rm (c)} The field $\Bbb F_p(\chi)$ satisfies $[\Bbb F_p(\chi):\Bbb F_p] = p^{2d}$, for an integer $d\ge1$, and $\chi(\vS)$ is contained in the subgroup of $\Bbb F_p(\chi)^\times$ of order $1{+}p^d$. 
\item 
Suppose that $\Bbb F_p(\chi)$ is symplectic. Any nonzero $\vS$-invariant alternating form on $\Bbb F_p(\chi)$ is $\vS$-anisotropic. Any two such forms are $\vS$-isometric. 
\item 
A finite $\Bbb F_p$-representation $U$ of $\vS$ provides a symplectic anisotropic rep\-resentation of $\vS$ if and only if there exist $\chi_j\in \Hom{}\vS{\bar{\Bbb F}_p^\times}$, $1\le j\le r$, such that 
\itemitem{\rm (a)} each $\Bbb F_p(\chi_j)$ is symplectic; 
\itemitem{\rm (b)} if $i\neq j$, then $\chi_i$ is not $\vO$-conjugate to $\chi_j$; 
\itemitem{\rm (c)} $U = \bigoplus_{1\le j\le r} \Bbb F_p(\chi_j)$. 
\endroster 
\endproclaim 
The proposition is taken from section 8.2 of \cite{1}. It may equally be viewed as an instance of the more general classification in \cite{9} Theorem 8.1, although some effort of translation would be required. 
\remark{Remark} 
If $\chi$ has order $a$ and satisfies the conditions in part (1) of the proposition, then 
\roster 
\item"\rm (a)" $a\ge3$ and 
\item"\rm (b)"  
the integer $d$ is the least for which $1{+}p^d$ is divisible by $a$. 
\endroster 
\endremark 
\subhead 
4.3 
\endsubhead 
In the same situation, we analyze the symplectic $\Bbb F_p$-representation of $\vS$ on $\vD = \Gal KT$. Following 4.2 Proposition (2), it is only the structure of the {\it linear\/} $\Bbb F_p\vS$-representation $\vD$ that need concern us. 
\par 
Recall that $J_{K/T}$ is the set of (upper) ramification jumps of $K/T$. For $j\in J_{K/T}$, define $W^j$ as in 3.2 Proposition. 
\proclaim{Proposition} 
For all $j\in J_{K/T}$, the $\Bbb F_p\vS$-space $W^j$ is irreducible. 
\endproclaim 
\demo{Proof} 
Let $k\in \Bbb Z$, $k\ge1$. The group $\vD^k$ is the image of the unit group $U^k_T$ under the Artin reciprocity map $T^\times \to \vD = \Gal KT$. This map is $\vS$-equivariant and $W^j$, $j\in J_{K/T}$, is so realized as a $\vS$-quotient of $U^j_T/U^{1+j}_T$. 
\par 
The natural action of $\vS = \Gal TF$ on $\frak p_T/\frak p_T^2$ is given by a faithful character $\theta:\vS \to \Bbbk_F^\times$. The natural action on $\frak p_T^j/\frak p_T^{1+j}$, $j\ge1$, is therefore implemented by $\theta^j$. The character $\theta^j$ induces an algebra homomorphism $\Bbb F_p\vS \to \Bbbk_F$, the image of which is necessarily a subfield of $\Bbbk_F$. The $\Bbb F_p\vS$-module $U^j_T/U^{1+j}_T \cong \frak p_T^j/\frak p_T^{1+j}$ is therefore {\it isotypic.} However, $W^j$ is anisotropic, so 4.2 Proposition (3) implies that $W^j$ is a direct sum of mutually inequivalent irreducible $\Bbb F_p\vS$-modules. It is therefore irreducible, as required. \qed 
\enddemo 
We underline some points made in the preceding proof. 
\proclaim{Corollary} 
Let $j\in J_{K/T}$. 
\roster 
\item 
The symplectic $\Bbb F_p$-representation $W^j$ of $\vS$ is equivalent to $\Bbb F_p(\theta^j)$. 
\item 
Let $e = e(T|F) = |\vS|$. An element $s\in \vS$ has a non-trivial fixed point in $W^j$ if and only if $s^{\gcd (e,j)} = 1$. 
\endroster 
\endproclaim 
\subhead 
4.4 
\endsubhead 
Twisting $\sigma$ with a character of $\Gal TF$ has no effect on the assertion to be proved. We therefore assume that the character $\det\sigma$ is trivial on $\vS$: this puts us in the situation of 2.2 Proposition.  
\par 
The orthogonal decomposition $\vD = \sum_{j\in J_{K/T}}  W^j$ implies a canonical realization of the $\Bbb F_p\vS$-module $W^j$ as a subspace of $\vD$: let $\wt W^j$ be its inverse image in $\wt\vD$. The construction outlined in 2.3, 2.4 gives a representation $\sigma^j$ of $\vS\wt W^j$ and a tensor decomposition 
$$ 
\sigma = \bigotimes_{j\in J_{K/T}} \sigma^j. 
$$ 
We may choose the factors $\sigma^j$ so that each character $\det\sigma^j\Mid \vS$ is trivial. As in 4.3 Corollary (2), $\vS$ acts on $W^j$ via its quotient of order $e_j = e/(e,j)$, and that action is faithful. 
\definition{Definition} 
Let $A$ be the set of positive divisors $a$ of $e$ of the form $e_j = e/(e,j)$, for some $j\in J_{K/T}$. For $a\in A$, set 
$$ 
\sigma_a = \bigotimes \Sb j\in J_{K/T}, \\ a = e_j \endSb \sigma^j. 
$$ 
\enddefinition 
Remark that $e = |\vS|$ is the lcm of the elements of $A$. Note also that a factor $\sigma_a$ may have several ramification jumps: this possibility is {\it not\/} excluded by 4.2 Proposition. 
\par 
We work in the ring $\Bbb Z\widehat \vS$ of virtual characters of $\vS$. The elements of $\Bbb Z\widehat \vS$ are thus the formal linear combinations 
$$ 
\bk c  = \sum_{\chi\in \widehat \vS} c_\chi\chi 
$$ 
in which the coefficients $c_\chi$ lie in $\Bbb Z$. Let $\Bbb N\widehat \vS$ be the ``order''  consisting of those $\bk c\in \Bbb Z\widehat \vS$ for which the coefficients $c_\chi$ are all non-negative. For $\bk a,\bk b \in \Bbb Z\widehat \vS$, we write $\bk a \ge \bk b$ when $\bk a{-} \bk b \in \Bbb N\widehat \vS$. We also use the relation $\ge$ to compare elements of $\Bbb Q\widehat\vS$ in the obvious way. 
\definition{Definition} 
Let $a \in A$. 
\roster 
\item Let $q_a$ be the least power of $p$ such that $1{+}q_a$ is divisible by $a$ and define $\ell(a)$ as the number of $j\in J_{K/T}$ for which $a = e_j$. 
\item Define the positive integer $k_a$ by 
$$ 
ak_a = \big(q_a^{\ell(a)}-(-1)^{\ell(a)}\big). 
$$ 
\item Let $\mu_a$ denote the trivial character of $\vS$ if $a$ is odd or $k_a$ is even. Otherwise, let $\mu_a \in \widehat \vS$ have order $2$. 
\item 
Let $\bs\rho_a\in \Bbb Z\widehat \vS$ be the sum of characters $\phi$ of $\vS$ such that $\phi^a = 1$, and define  
$$ 
R_a = k_a\bs\rho_a + (-1)^{\ell(a)}\mu_a \in \Bbb Z\widehat \vS. 
\tag 4.4.1 
$$ 
\endroster 
\enddefinition 
\proclaim{Proposition} 
If $a\in A$, then 
$$ 
\align 
\sigma_a\Mid \vS &= R_a \tag 4.4.2 \\ 
\intertext{and, moreover,} 
\sigma\Mid\vS &= \prod_{a\in A} R_a, \tag 4.4.3 
\endalign 
$$ 
the product being taken in $\Bbb Z\widehat\vS$. 
\endproclaim 
\demo{Proof} 
This follows directly from 2.2 Proposition and (2.4.1). \qed 
\enddemo 
\subhead 
4.5 
\endsubhead 
We treat a special case of 4.1 Theorem, working directly from (4.4.2). 
\proclaim{Proposition} 
If the set $A$ has exactly one element, then $m(\sigma) \le \sw(\sigma)$. 
\endproclaim 
\demo{Proof} 
The lcm of the elements of $A$ is $e$, so $A = \{e\}$. Moreover, $q_e^{\ell(e)} = p^r = \dim\sigma$. 
\par 
Suppose first that $\ell(e)$ is odd, so that $m(\sigma) = k_e = (p^r{+}1)/e$. We have to show that $p^r{+}1 \le e\,\sw(\sigma)$. By 3.3 Theorem, $e\,\sw(\sigma) \ge (p^r{+}1)j_\infty$, where $j_\infty$ is the largest element of $J_{K/T}$. As $j_\infty$ is a positive integer ({\it cf\.} 3.2), so 
$$ 
e\,\sw(\sigma) \ge (p^r{+}1)j_\infty \ge (p^r{+}1), 
$$ 
as required. Suppose, on the other hand, that $\ell(e)$ is even. In this case, $e$ divides $p^r{-}1$, so $e\le p^r{-}1$ and 
$$ 
em(\sigma) =p^r{-}1{+}e \le 2(p^r{-}1).  
$$ 
On the other hand, as $\ell(e)$ is even so $j_\infty \ge2$. Therefore 
$$ 
em(\sigma) \le 2(p^r{-}1) < (p^r{+}1)j_\infty \le e\sw(\sigma). 
$$ 
This completes the proof. \qed 
\enddemo 
We reflect briefly on the proof of this proposition. 
\proclaim{Corollary of proof} 
In the situation of the proposition, if $m(\sigma) = \sw(\sigma)$ then $J_{K/T} = \{1\}$. 
\endproclaim 
\demo{Proof} 
If $\ell(e)$ is even then, as we have just seen, $m(\sigma) < \sw(\sigma)$, so suppose $\ell(e)$ is odd. If $\ell(e) \neq 1$, then $j_\infty \ge 3$ and 
$$ 
em(\sigma) = 1{+}p^r < (1{+}p^r)j_\infty \le e\,\sw(\sigma). 
$$ 
So, we assume $\ell(e) = 1$ and $J_{K/T} = \{j_\infty\}$. In this case, if $j_\infty > 1$ then $m(\sigma) < \sw(\sigma)$. \qed 
\enddemo 
\subhead 
4.6 
\endsubhead 
Assume now that $A$ has at least two elements. We can make some simplifying approximations. The expressions (4.4.1), (4.4.2) imply 
$$ 
\sigma_a\Mid \vS = R_a \le q_a^{\ell(a)}\,d_a\,\frac{\bs\rho_a}a 
\tag 4.6.1 
$$ 
where 
$$ 
d_a = \left\{\,\alignedat3 &1+q_a^{-\ell(a)} &\quad &\text{if $\ell(a)$ is odd,} \\&1 + (a{-}1)q_a^{-\ell(a)} &\quad &\text{if $\ell(a)$ is even.} 
\endalignedat \right. 
\tag 4.6.2 
$$ 
If $a,b\in A$, then 
$$ 
\frac{\bs\rho_a}a\,\frac{\bs\rho_b}b = \frac{\bs\rho_c}c, 
$$ 
where $c$ is the lcm of $a$ and $b$. So, taking the product over $a\in A$, we get 
$$ 
\sigma\Mid \vS \le \dim\sigma \prod_{a\in A} d_a\,\frac{\bs\rho_e}e, 
$$ 
whence 
$$ 
m(\sigma) \le \frac{\dim\sigma}e\,\prod_{a\in A} d_a. 
\tag 4.6.3 
$$ 
So, we are reduced to proving: 
\proclaim{Proposition} 
If\/ $|A| \ge2$, then 
$$ 
\prod_{a\in A} d_a \le \frac{e\,\sw(\sigma)}{\dim\sigma}. 
\tag 4.6.4 
$$ 
\endproclaim 
\demo{Proof} 
Let $j_\infty$ be the largest element of $J_{K/T}$. Let $\vX$ be a $J$-split Lagrangian subspace of the symplectic space $\vD$ as in 3.3. Let $L$ be the fixed field of $\vX$ and recall that $J_{K/T}$ is equal to the set $J_{L/T}$ of jumps of $L/T$. 
\proclaim{Lemma 1} 
Suppose that $A$ has at least two elements. If $a\in A$, then $a\ge 3$ and $d_a \le a/(a{-}1)$. 
\endproclaim 
\demo{Proof} 
For the first assertion, see 4.2 Remark (a). For the second, suppose first that $\ell(a)$ is even, whence $d_a = 1{+}(a{-}1)q_a^{-\ell(a)}$. By definition, $a$ divides $1{+}q_a$ whence $a{-}1\le q_a$. Therefore 
$$ 
d_a \le 1+q_a^{-1} \le 1+(a{-}1)^{-1} = a/(a{-}1). 
$$ 
The case of $\ell(a)$ odd is similar but even easier. \qed 
\enddemo 
It follows that 
$$ 
\prod_{a\in A} d_a \le \prod_{a\in A} \frac a{a{-}1}. 
$$ 
Let $t\ge2$ be the number of elements of $A$. The largest element of $A$ is therefore, at least, $t{+}2$ and all elements are $\ge3$. Consequently, 
$$ 
\prod_{a\in A} d_a \le \frac{3}{2}\, \frac{4}{3}\,\cdots \, \frac{t{+}2}{t{+}1}\,\le\,\frac{t{+}2}2 \,\le\, t. 
$$ 
\proclaim{Lemma 2} 
If $h$ denotes the number of jumps of $L/T$, then 
$$ 
h < p^{-r}\psi_{L/T}(j_\infty) + j_\infty. 
$$ 
\endproclaim 
\demo{Proof} 
The jumps of the abelian extension $L/T$ are positive integers, so $h\le j_\infty$ and $h-j_\infty \le 0 < p^{-r}\psi_{L/T}(j_\infty)$, as required. \qed 
\enddemo 
Assembling these relations and applying (3.3.1), we get 
$$ 
\prod_{a\in A} d_a \le t \le h < e\,\sw(\sigma)/p^r, 
$$ 
as required to complete the proof of the proposition. \qed
\enddemo 
This finishes the proof of 4.1 Theorem. \qed 
\head 
5. An estimate of the different 
\endhead 
Preliminary to the proof of the general case of the main theorem, we make an estimate of the wild exponent $w_{K/F}$ of a class of finite extensions $K/F$. It is not remotely sharp (see 5.2 Example) but is adequate for our purposes. 
\subhead 
5.1 
\endsubhead 
Let $K/F$ be a finite separable extension, with $K\subset \bar F$. The {\it wild exponent\/} $w_{K/F}$ of $K/F$  is 
$$ 
\aligned 
w_{K/F} &= d_{K/F} +1 - e(K|F) \\ 
&= \sw\big(\Ind_{K/F}\,1_K\big), 
\endaligned 
\tag 5.1.1 
$$ 
where $d_{K/F}$ is the exponent of the different of $K/F$ and $1_K$ denotes the trivial character of $\scr W_K$. 
\subhead 
5.2 
\endsubhead 
Let $E/F$ be a finite, totally ramified, Galois extension. Set $\Gal EF = \vG$ and let $\vD$ be the wild inertia subgroup of $\vG$. As in 1.1, $\vD$ is the unique $p$-Sylow subgroup of $\vG$ and admits a complement $\vS$ in $\vG$. In particular, $\vS$ is cyclic of order prime to $p$. 
\proclaim{Proposition} 
Let $\vF$ be a subgroup of $\vG$, such that the index $(\vG:\vF)$ is a power of $p$. If $K$ is the fixed field $E^\vF$ of $\vF$ in $E$, then 
$$ 
w_{K/F} \ge \big|\vF\backslash \vG/\vS\big| - 1. 
$$ 
\endproclaim 
\demo{Proof} 
Recall that any two choices of the complement $\vS$ are conjugate in $\vG$. The assertion is therefore independent of the choice of $\vS$. 
\par 
If $\vF = \vG$ there is nothing to prove, so we assume otherwise. 
\proclaim{Lemma 1} 
Let $\vX$ be a normal subgroup of $\vG$ such that $\vX \subset \vF$ and let $f:\vG\to \vG/\vX$ be the quotient map. 
\roster 
\item The group $f(\vS)$ is a complement of $f(\vD)$ in $\vG/\vX$. 
\item 
The map $f$ induces a $\vS$-equivariant bijection $\vF\backslash \vG \to f(\vF)\backslash f(\vG)$, and hence a bijection $\vF\backslash \vG/\vS \to f(\vF)\backslash f(\vG)/f(\vS)$. 
\endroster 
\endproclaim 
\demo{Proof} 
Straightforward. \qed 
\enddemo 
Continue with $\vX$ as in the lemma. If we replace $E$ by $E^\vX$, the extension $K/F$ is unchanged. The effect of the lemma is to show that, if the proposition holds for the configuration $F\subset K\subset E^\vX$, then it holds for $F\subset K\subset E$. We may choose $\vX$ so that $E^\vX/F$ is a normal closure of $K/F$. It is therefore enough to prove the proposition under the assumption that $E/F$ {\it is a normal closure of\/} $K/F$. We henceforward assume this to be the case. 
\proclaim{Lemma 2}  
Let $\vT$ be the smallest non-trivial ramification subgroup of $\vG$. The group $\vT$ is elementary abelian and central in $\vD$. It is not contained in $\vF$. 
\endproclaim 
\demo{Proof} 
The first assertions are given by \cite{12} IV Prop\. 7 and Prop\. 10. If $\vT$ were contained in $\vF$ then $E^\vT/F$ would be a normal extension containing $K/F$ and such that $[E^\vT:F] < [E:F]$. Since $E/F$ is a normal closure of $K/F$, this is impossible. \qed 
\enddemo \enddemo 
Suppose for the moment that $\vG = \vF\vT$ or, equivalently, that $\vD = (\vF\cap \vD)\vT$. As $\vT$ is central in $\vD$, so $\vF\cap \vD$ is normal in $\vD$ and $\vD/\vF\cap \vD$ is abelian. Let $1_\vF$ denote the trivial character of $\vF$, and similarly for other groups. The Mackey formula gives the relations 
$$ 
\align 
\Ind_\vF^\vG\,1_\vF\Mid \vD &= \Ind_{\vF\cap \vD}^\vD\,1_{\vF\cap \vD}, \tag 5.2.1 \\
\Ind_\vF^\vG\,1_\vF\Mid \vT &= \Ind_{\vF\cap \vT}^\vT\,1_{\vF\cap \vT}. \tag 5.2.2 
\endalign 
$$ 
The restriction (5.2.2) is the direct sum of all characters $\chi$ of $\vF{\cap} \vT\backslash \vT$. Any such character $\chi$ extends uniquely to a character $\chi_\vD$ of $\vD$ trivial on $\vF{\cap} \vD$: one puts $\chi_\vD(hr) = \chi(r)$, $h\in \vF{\cap} \vD$, $r\in \vT$. Consequently, 
$$ 
\Ind_{\vF\cap \vD}^\vD \,1_{\vF\cap \vD} = \sum_{\chi \in (\vF\cap \vT\backslash \vT)\sphat}\chi_\vD. 
$$ 
If $\vG_\chi$ denotes the $\vG$-centralizer of $\chi_\vD \in (\vD/\vF{\cap} \vD)\sphat$, then $\vG_\chi = \vS_\chi \vD$, where $\vS_\chi$ is the $\vS$-centralizer of $\chi_\vD$ (or, equivalently, of $\chi$). Consequently, there is a unique character $\chi_\vS$ of $\vG_\chi$ that extends $\chi_\vD$ and is trivial on $\vS_\chi$. Therefore 
$$ 
\Ind_\vF^\vG\,1_\vF = \sum_{\chi\in \vS\backslash (\vF\cap \vT\backslash \vT)\sphat}\ \sum_{\eta\in (\vG_\chi/\vD)\sphat} \Ind_{\vG_\chi}^\vG\, \eta\chi_\vS. 
$$ 
We calculate the contribution of each term here to the exponent $\sw(\Ind_\vF^\vG\,1_\vF) = w_{K/F}$. 
\par 
If $\chi$ is trivial, then $\vG_\chi = \vG$ and we get a contribution of $0$. Otherwise, $\Ind_{\vG_\chi}^\vG\, \eta\chi_\vS$ has Swan exponent at least $1$, whence 
$$ 
w_{K/F} \ge \sum \Sb \chi\in \vS\backslash (\vF\cap \vT\backslash \vT)\sphat \\ \chi\neq 1 \endSb (\vG_\chi:\vD) . 
$$ 
However, $\big|\vS\backslash (\vF\cap \vT\backslash \vT)\sphat{\,\,}\big| = \big|\vF\backslash \vG/\vS\big|$, and so $w_{K/F} \ge \big|\vF\backslash \vG/\vS\big|-1$ in this case. 
\example{Example}
Remark here that, once the trivial character $\chi$ is excluded, all groups $\vG_\chi$ are the same: they depend only on the denominator of $j$, where $\vT = \vG^j \neq \vG^{j+\eps}$, $\eps >0$. All characters $\eta\chi_\vS$ have the same slope, namely $j$. The index $(\vG:\vG_\chi)$, for $\chi\neq 1$, is the g.c.d\. of $|\vS|$ and the denominator of $j$. So, for $\chi\neq 1$, the inner sum has Swan exponent $j(\vG:\vG_\chi) (\vG_\chi:\vD) = j|\vS|$. Therefore 
$$ 
\sw(\Ind_\vF^\vG\,1_\vF) = w_{K/F} =  j|\vS|\,\big(\big|\vF\backslash \vG/\vS\big| -1\big). 
\tag 5.2.3 
$$ 
\endexample  
We return to the proof of 5.2 Proposition, assuming now that $\vF\vT \neq \vG$. Since the index $(\vG:\vF)$ is a power of $p$, the group $\vF$ contains a conjugate of $\vS$. Following the remark at the beginning of the proof, we may assume that $\vS\subset \vF$. 
\par 
Let $L = E^{\vF\vT}$. The first case above gives 
$$ 
w_{K/L} \ge \big|\vF\backslash \vF\vT/\vS\big| - 1. 
$$ 
By induction on $[K{:}F] = (\vG:\vF)$, we likewise have 
$$ 
w_{L/F} \ge \big|f(\vF)\backslash f(\vG)/f(\vS)\big|-1, 
$$ 
where $f:\vG \to \vG/\vT$ is the quotient map. On the other hand, 
$$ 
w_{K/F} = w_{K/L} + [K{:}L]\,w_{L/F}, 
$$ 
so 
$$ 
w_{K/F} \ge \big|\vF\backslash \vF\vT/\vS\big|-1 + [K{:}L]\big(\big|f(\vF)\backslash f(\vG)/f(\vS)\big| - 1\big). 
$$ 
Under the canonical surjection $\bar f: \vF\backslash \vG/S \to f(\vF)\backslash f(\vG)/f(\vS)$ induced by the quotient map $f:\vG \to \vG/\vT$, the fibre of the trivial coset $f(\vF) = f(\vF)f(\vS)$ is precisely $\vF\backslash \vF\vT/\vS$. On the other hand, let $x = f(g) \notin f(\vF)$. The fibre, under $\bar f$, of $f(\vF)xf(\vS)$ is contained in $\vF g\vS$. This comprises at most $[K{:}L]$ double cosets $\vF g\vS$, whence the result follows. \qed 
\head 
6. Proof of the main theorem 
\endhead 
We prove the Tame Multiplicity Theorem in the general case. Let $\sigma$ be an irreducible representation of $\scr W_F$ that is {\it not tamely ramified:\/} see 1.1 Remark (1). Since the assertion of the theorem is unaffected by tensoring $\sigma$ with an unramified character of $\scr W_F$, we may treat $\sigma$ as a representation of $\vG = \Gal EF$, where $E/F$ is finite. Let $\vG_0$, $\vG_1$ be respectively the inertia and the wild inertia subgroups of $\vG$, and similarly for other finite Galois groups. 
\subhead 
6.1 
\endsubhead 
Let $\sigma$ be an irreducible representation of $\vG = \Gal EF$, with $\sw(\sigma) > 0$. Let $\vS$ be a complement of $\vG_1$ in $\vG_0$. 
\proclaim{Proposition} 
If $\sigma$ is absolutely ramified, that is, if $E/F$ is totally ramified, then $m(\sigma) \le \sw(\sigma)$. 
\endproclaim 
\demo{Proof} 
If $\sigma$ is primitive or of dimension one, the result holds by 4.1 Theorem or 1.1 Remark (2) respectively. We therefore suppose otherwise: there is a proper subgroup $\vD$ of $\vG$ and an irreducible representation $\tau$ of $\vD$ such that $\sigma  = \Ind_\vD^\vG\,\tau$. The representation $\tau$ is absolutely ramified and, by induction on dimension, we may assume that $m(\tau) \le \sw(\tau)$. 
\par 
Suppose first that $\vD$ may be chosen to contain $\vG_1$. Thus $\vG = \vG_0 = \vS\vD$, and $\vD$ is a normal subgroup of $\vG$. The Mackey formula gives 
$$ 
\sigma\Mid \vS = \Ind_{\vS\cap \vD}^\vS\,\tau\Mid \vS\cap \vD = \tau\Mid\vS, 
$$ 
whence $m(\sigma) = m(\tau)$. As $E^\vD/F$ is tamely ramified, so $\sw(\sigma) = \sw(\tau)$ and we are done in this case. 
\par 
We therefore assume that $\sigma$ cannot be induced from a proper subgroup of $\vG = \vG_0$ that contains $\vG_1$. Since $\vG/\vG_1$ is cyclic, the restriction $\sigma\Mid \vG_1$ is irreducible. In particular, $\dim\sigma$ is a power of $p$. It follows that, if $\sigma$ is induced from a representation $\tau$ of a proper subgroup $\vD$ of $\vG$, then $(\vG{:}\vD)$ is a power of $p$ and, if $K = E^\vD$, the extension $K/F$ is totally wildly ramified. We have $m(\tau) \le \sw(\tau)$, and 
$$ 
\sw(\sigma) = \sw(\tau)+w_{K/F}\,\dim\tau. 
\tag 6.1.1 
$$ 
We apply 5.2 Proposition. We adjust our choice of $\vS$, via conjugation by an element of $\vG_1$, to achieve $\vS\subset \vD$. Let $\chi$ be a character of $\vS$. In the Mackey expansion 
$$ 
\sigma\Mid \vS = \sum_{g\in \vD\backslash \vG/\vS} \Ind_{g^{-1}\vD g\cap \vS}^\vS\,\big(\tau^g\Mid g^{-1}\vD g\cap \vS\big). 
$$ 
the trivial double coset gives the term $\tau\Mid \vS$, in which $\chi$ occurs with multiplicity at most $m(\tau)$. The contribution from a non-trivial double coset contains $\chi$ with multiplicity at most $\dim\tau$ so, overall, 
$$ 
m(\sigma) \le m(\tau) + \big(|\vD\backslash\vG/\vS|-1\big)\,\dim\tau. 
\tag 6.1.2 
$$ 
Comparing (6.1.1) with (6.1.2), the proposition of 5.2 implies 
$$ 
m(\sigma) \le \sw(\tau) + w_{K/F}\,\dim\tau = \sw(\sigma), 
$$ 
as required. \qed 
\enddemo 
We return to a more general situation. 
\proclaim{Corollary} 
Let $E/F$ be a finite Galois extension and let $\sigma$ be an irreducible representation of $\vG = \Gal EF$, with $\sw(\sigma) >0$. If $\sigma\Mid \vG_0$ is irreducible, then $m(\sigma)\le \sw(\sigma)$. 
\endproclaim 
\demo{Proof} 
The representation $\sigma_0 = \sigma\Mid\vG_0$ is irreducible and absolutely ramified. The proposition gives $m(\sigma_0) \le \sw(\sigma_0)$. However, since $\vS\subset \vG_0$, we have $m(\sigma_0) = m(\sigma)$. On the other hand, $\sw(\sigma_0) = \sw(\sigma)$, since $E^{\vG_0}/F$ is unramified. \qed 
\enddemo 
\example{Example} 
Example 2 of \cite{3} 8.5 is interesting in this context. Suppose that $p=2$ and that $F$ contains a primitive cube root of unity. The construction in \cite{3} yields a primitive representation $\sigma$ of dimension $8$, with $\sw(\sigma) = 3$ and a unique ramification jump. (In the notation of 3.3 Theorem, this jump is $j_\infty$ and it has value $1$.) If $\roman{Ker}\,\bar\sigma = \scr W_K$, and $T/F$ is the maximal tame sub-extension of $K/F$, then $[T{:}F] = 9$ and $e(T|F) = 3$. In particular, $\sigma$ is not absolutely ramified. If $T_0/F$ is the maximal unramified sub-extension of $T/F$, the restriction $\sigma\Mid \scr W_{T_0}$ is irreducible but not primitive. A simple counting argument gives $m(\sigma) =  3 = \sw(\sigma)$. 
\endexample 
\subhead 
6.2 
\endsubhead 
We complete the proof of the Tame Multiplicity Theorem. Let $\sigma$ be an irreducible representation of the finite group $\vG = \Gal EF$ with $\sw(\sigma) > 0$. Let $\vS$ be a complement of $\vG_1$ in $\vG_0$. If $\sigma \Mid \vG_0$ is irreducible, the theorem is 6.1 Corollary. We therefore assume otherwise, so there exist a proper subgroup $\vD$ of $\vG$ containing $\vG_0$ and an irreducible representation $\tau$ of $\vD$ that induces $\sigma$. We choose $\vD$ minimal with respect to this property, so that $\tau \Mid \vD_0$ is irreducible. By 6.1 Corollary, $m(\tau) \le \sw(\tau)$ while 
$$ 
\sw(\sigma) = (\vG{:}\vD)\,\sw(\tau). 
\tag 6.2.1 
$$ 
As $\vD_0 = \vG_0$ and $\vD_1 = \vG_1$, so $\vS$ is also a complement of $\vD_1$ in $\vD_0$. Applying the standard Mackey formula, we get 
$$ 
\sigma\Mid \vS = \sum_{g\in \vD\backslash \vG/\vS} \Ind_{g^{-1}\vD g\cap \vS}^\vS\,\big(\tau^g\Mid g^{-1}\vD g\cap \vS\big). 
$$ 
We have $\vG_0 = \vS \vG_1 \subset \vD$, while any $\vG$-conjugate of $\vS$ is contained in $\vD$. The canonical map $\vD\backslash \vG \to \vD\backslash \vG/\vS$ is therefore bijective. Consider the expression 
$$ 
\sigma\Mid \vS = \sum_{g\in \vD\backslash \vG} \tau^g\Mid \vS. 
$$ 
If $\chi$ is a character of $\vS$, the multiplicity of $\chi$ in $\tau^g$ is that of $\chi^{g^{-1}}$ in $\tau$, whence at most $m(\tau)$. We conclude that $m(\sigma) \le (\vG{:}\vD)\,m(\tau)$. Since $m(\tau) \le \sw(\tau)$, the desired relation $m(\sigma) \le \sw(\sigma)$ follows from (6.2.1). \qed 
 

\Refs 
\ref\no1 
\by C.J. Bushnell and A. Fr\"ohlich 
\book Gauss sums and $p$-adic division algebras 
\bookinfo Lecture Notes in Math. \vol 987 \publ Springer \publaddr Berlin-Heidelberg New York  \yr 1983 
\endref 
\ref\no2 
\by C.J. Bushnell and G. Henniart 
\paper 
Langlands parameters for epipelagic representations of  $\roman{GL}_n$ 
\jour Math. Annalen \vol 358 \yr 2014 \pages 433--463 
\endref 
\ref\no3
\bysame 
\paper Higher ramification and the local Langlands correspondence \jour Annals of Mathematics (2) \vol 185 \issue 3 \yr 2017 \pages 919--955 
\endref 
\ref\no4 
\bysame 
\paper Local Langlands correspondence and ramification for Carayol representations \jour arXiv: 1611.09258v3. To appear in Compositio Math  
\endref 
\ref\no5 
\by G. Glauberman 
\paper Correspondences of characters for relatively prime operator groups 
\jour Canadian J. Math. \vol 20 \yr 1968 \pages 1465--1488 
\endref 
\ref\no6 
\by D. Gorenstein 
\book Finite groups 
\publ AMS Chelsea Publishing \publaddr Providence RI \yr 2012 
\endref 
\ref\no7 
\by G. Henniart 
\paper Repr\'esentations du groupe de Weil d'un corps local  \jour L'Ens. Math. S\'er. II \vol 26 \yr 1980 \pages 155-172 
\endref 
\ref\no8 
\by I.M. Isaacs 
\book Character theory of finite groups 
\publ AMS Chelsea Publishing \publaddr Providence RI \yr 2011  
\endref 
\ref\no9 
\by H. Koch 
\paper Classification of the primitive representations of the Galois group of local fields 
\jour Invent. Math. \vol 40 \yr 1977 \pages 195--216 
\endref 
\ref\no10 
\by M. Reeder 
\paper Adjoint Swan conductors I: the essentially tame case 
\jour Int. Math. Res. Not. IMRN \yr  2017  \pages doi: 10.1093/imrn/rnw301 
 \endref 
\ref\no11 
\by J.F. Rigby 
\paper Primitive linear groups containing a normal nilpotent subgroup larger than the centre of the group 
\jour J. London Math. Soc. \vol 35 \yr 1960 \pages 389--400 
\endref 
\ref\no12  
\by J-P. Serre 
\book Corps locaux \publ Hermann \publaddr Paris \yr 1968 
\endref 
\endRefs
\enddocument